# Novel numerical method for simultaneous design and control optimization of seasonal thermal energy storage systems

Wonsun Song [a] [✉] [*], Jakob Harzer [a] [✉], Christopher Jung [b] [✉], Leon Sander [b] [✉], Moritz Diehl [a,c]

[a] *Department of Microsystems Engineering (IMTEK), University of Freiburg, 79110, Freiburg, Germany*
[b] *Department of Environmental Meteorology, University of Freiburg, 79085, Freiburg, Germany*
[c] *Department of Mathematics, University of Freiburg, 79104, Freiburg, Germany*

## ARTICLE INFO

## ABSTRACT

The transition to a carbon-neutral energy system requires widespread deployment of renewable energy sources and economically feasible energy storage solutions. This study presents a comprehensive optimization framework that jointly addresses the design and control of a nonlinear energy system supplying both heat and electricity to the Dietenbach district in Freiburg, Germany. The proposed system integrates solar and wind power with battery storage and seasonal thermal energy storage coupled via a heat pump, enhancing self-sufficiency and mitigating seasonal supply–demand mismatches. A multi-node lumped-parameter model captures heat transfer within the pit thermal energy storage, forming the basis of a periodic optimal control problem solved numerically. An averaging method reduces computation time by 80.5% while preserving fidelity for year-long optimization. A case study shows a projected total yearly energy cost of 5.93 €/m² for combined heat and electricity, which is 73% lower than the German average. This study underscores the feasibility of designing economically viable, autonomous energy communities in real-world scenarios and provides an efficient, robust optimization framework for designing system components and operational control strategies.

## 1. Introduction

Achieving carbon-neutral district energy systems calls for both efficient technology choices and intelligent operational strategies. While the transition to renewables has made significant strides — Germany, for instance, now sources over 50% of its electricity from renewables — emissions from the heating sector remain stubbornly high [1]. The reliance of this sector on fossil fuels, the apparent seasonal mismatch between energy supply and demand, and the inherent difficulty of transporting heat over long distances underscore the need for regionally tailored large-scale thermal energy storage solutions.

Pit Thermal Energy Storage (PTES) provides a promising route to long-term heat storage, offering high charging and discharging rates, flexible capacity, and cost-effectiveness—qualities that make it particularly suitable for cross-sectoral energy systems [2]. Meanwhile, recent large-scale heat pump applications demonstrate the feasibility of using a heat pump as a primary heat source, overcoming inherent drops in the coefficient of performance (COP) [3]. Integrating these technologies into a single system can significantly reduce fossil fuel use, capitalize on sector coupling, and address both short-term fluctuations and seasonal imbalances.

However, the design and control of PTES with source-side heat pumps remain unexplored in the literature, in part because capturing their transient, temperature-dependent dynamics at a year-long scale is computationally challenging. Many existing approaches avoid these nonlinearities (resulting from the joint control and design and from the heat flow dynamics) or rely on scenario-based simulations, missing opportunities for fully optimized, derivative-based solutions. Consequently, a framework that can simultaneously capture system design, nonlinear thermodynamics, and a long-horizon operational strategy is still lacking.

This work addresses these gaps by proposing an integrated energy system that combines photovoltaics, wind energy, PTES, and battery storage with a large-scale heat pump at its source side. The objective is

---

| Abbreviations | |
|---|---|
| OCP | Optimal Control Problem |
| LP | Linear Programming |
| NLP | Nonlinear Programming |
| MINLP | Mixed Integer Nonlinear Programming |
| MILP | Mixed Integer Linear Programming |
| STES | Seasonal Thermal Energy Storage |
| PTES | Pit Thermal Energy Storage |
| BTES | Borehole Thermal Energy Storage |
| PV | Photovoltaics |
| HP | Heat Pump |
| COP | Coefficient of Performance |
| SOC | State of Charge |
| ANI | Annuity of the Investment |
| ANF | Annuity Factor |
| CAPEX | Capital Expenditures |
| OPEX | Operational Expenditures |
| RMSLE | Root Mean Squared Logarithmic Error |
| RMSE | Root Mean Squared Error |
| NRMSE | Normalized Root Mean Squared Error |
| GHI | Global Horizontal Irradiance |
| GSI | Global Solar Irradiance |
| PR | Performance Ratio |
| SH | Space Heating |
| DHW | Domestic Hot water |

| Optimization symbol | |
|---|---|
| $x$ | State (–) |
| $u$ | Control (–) |
| $s$ | Scaling factor (–) |
| $\theta$ | Parameter vector (–) |
| $h$ | Step size (–) |
| $t$ | Time (s) |
| $N_c$, $N_f$ | Number of Coarse, Fine discretization Points (–) |
| $n$ | Number of variables (–) |

| Heat transfer symbol | |
|---|---|
| $\lambda$ | Thermal conductivity (W/(m K)) |
| $\dot{Q}$ | Heat Transfer Rate (W) |
| $P$ | Power (W) |
| $U$ | Overall heat transfer coefficient (W/(m² K)) |
| $R$ | Thermal resistance (K/W) |
| $C_s, C_g$ | Thermal capacitance (J/K) |
| $c_p$ | Specific heat capacity (J/(kg K)) |
| $T$ | Temperature (K) |
| $\dot{T}$ | Rate of temperature change (K/s) |
| $\dot{m}$ | Mass flow rate (kg/s) |
| $M$ | Number of storage layer (–) |
| $N$ | Number of ground layer (–) |
| $V$ | Volume (m³) |
| $A$ | Area (m²) |
| $z$ | Storage height (m) |
| $d$ | Ground distance (m) |
| $\eta$ | Efficiency (–) |
| $\rho$ | Density (kg/m³) |

| Miscellaneous symbol | |
|---|---|
| $I$ | Investment Cost (€) |
| $c$ | Price (€/kWh) |
| $n$ | Years (–) |
| $r$ | Interest rate (%) |
| $C$ | Installed capacity (W or Wh) |

| Indices | |
|---|---|
| f | Fine |
| c | Coarse |
| re | Renewable Energy |
| b | Battery |
| ch | Charging |
| dis | Discharging |

| s | Storage |
|---|---|
| g | Ground |
| bc | Boundary Condition |
| el | Electricity |
| hh | Households |
| sup | Supply |
| ret | Return |
| eff | Effective |
| top | Top layer |
| amb | Ambient |
| q | Cross-sectional |
| surf | Surface |

to meet a community's electricity and heating demands with minimal dependence on the electricity grid. This raises two key questions:

1. For a given demand profile over a year, how should the energy flows be controlled to minimize the dependence on the costly electricity grid?
2. How should the system components be sized to minimize the running and investment costs?

The answers to these control and design questions are strongly coupled and must therefore be addressed simultaneously. To this end, a single, large nonlinear optimization problem is constructed:

$$\underset{\substack{\text{(control strategy)}\\\text{(design parameters)}}}{\text{minimize}} \quad \text{(total cost)} \tag{1a}$$

$$\text{subject to} \quad \text{(year-periodic system simulation),} \tag{1b}$$

$$\text{(operational constraints),} \tag{1c}$$

$$\text{(demand satisfaction)} \tag{1d}$$

and then solved numerically using a full year of data on energy demand, renewable electricity supply, and ambient temperature.

### 1.1. Related work

There exists a variety of optimization frameworks which support multi-energy system modeling. Energy Hub Design Optimization (EHDO) targets hub configuration and operation planning [4]. Energy-PLAN assesses renewable energy integration and flexibility at regional and national scales [5]. Balmorel focuses on the electricity and district heating sectors with a bottom-up approach [6]. OSeMOSYS (Open Source energy MOdelling SYStem) is widely used for long-term policy and investment scenario analysis [7]. oemof (Open Energy Modeling Framework) provides modular, component-based energy-system models [8]. PyPSA (Python for Power System Analysis) focuses on power network studies, including multi-carrier flows and grid expansion [9].

To ensure computational tractability, these frameworks typically rely on linear or mixed-integer linear programming (LP or MILP) formulations, such as simplified discrete representations of energy storage that neglect temperature and mass flow dynamics. While effective for high-level planning, such simplifications limit their capability to investigate scenarios and control strategies involving transient ground interaction and thermal stratification inherent in PTES and heating systems.

For detailed physical modeling, the transient behavior of PTES is usually described by partial differential equations solved with commercial software. COMSOL Multiphysics applies the finite element method and delivers high spatial fidelity [10], but is ill-suited for optimization of PTES within renewable energy systems. TRNSYS adopts a coarse-resolution model that eases system integration [11], yet remains focused on simulation rather than optimization. Both tools are proprietary, restricting source code access and limiting opportunities for customization, and seamless integration into a derivative-based optimization framework is nontrivial.

Recent contributions have advanced the state of PTES modeling. Dahash et al. developed an efficient finite element model and presented guidelines for the numerical simulation of PTES [12]. Ochs et al. compared different modeling approaches, quantifying the trade-off between accuracy and computational cost [13]. In a comprehensive review, Xiang et al. covered key aspects of PTES, including geometry, material selection, numerical methods, and demonstration projects [14]. Although these studies provide valuable insights into PTES system, they also highlight a persistent gap: the lack of open-source, computationally tractable models specifically designed for system-level optimization.

This gap forces an apparent dilemma in optimization between the modeling horizon and physical detail. On the one hand, studies using nonlinear or mixed-integer nonlinear programming (NLP or MINLP) to capture detailed physics are confined to short operational horizons. For example, Lu et al. solved a day-ahead MINLP to schedule building energy systems with thermal storage [15]. Jansen et al. introduced a MINLP model-predictive controller for a district heating network, demonstrating feasibility only for short-term control windows [16]. Extending such nonlinear models to a full-year horizon can require several days of computation for a single run, as shown by Walden et al. in a power-to-heat case [17].

On the other hand, achieving a year-long optimization horizon for seasonal storage has so far been possible only with highly simplified formulations, mainly focused on borehole thermal energy storage (BTES). Zhou et al. for example, developed a MILP for an integrated electricity and heat system, which ignored the mass flow and temperature dynamics [18]. In a similar vein, Yu et al. employed an LP model for co-optimizing BTES and power-to-gas units, which also omitted any dependence on temperature or mass flow [19]. Addressing some of these limitations, Fiorentini et al. proposed a mixed-integer, bilinear (non-convex) model [20]. However, this approach still relied on a highly aggregated representation, tracking only one lumped BTES temperature, one thermal capacitance, and a single ground-temperature state, which compromised the physical accuracy of the heat transfer dynamics.

To ease the computational burden of long-horizon optimization, many studies reduce the problem's temporal resolution. Kotzur et al. linked the state of charge of storage in aggregated typical periods, preserving seasonal storage behavior [21]. Gabrielli et al. employed a single, highly aggregated storage variable for the entire year [22]. Blanke et al. reviewed both approaches and extended Kotzur's model by merging consecutive identical typical periods [23]. Schütz et al. compared clustering algorithms for generating typical demand days and identified k-medoids as the most accurate choice for energy system [24]. Building on this body of work, this paper proposes a novel averaging-based method that is directly integrated into the optimization process, improving both tractability and model fidelity.

To the best of the authors' knowledge, no open-source optimization framework captures the nonlinear behavior of heating networks and thermal energy storage, while simultaneously integrating system design and operational control in a year-long horizon optimization for multi-energy systems—all while maintaining computational efficiency.

### 1.2. Contribution

This study addresses the identified research gaps by formulating an optimal control problem (OCP) that accurately captures key system dynamics while ensuring computational feasibility. The primary contributions are:

1. Advanced, yet optimization-friendly thermal energy storage modeling: A stratified underground thermal energy storage model is developed based on a multi-node lumped-parameter method. This approach captures essential dynamics, such as temperature stratification and mass flow, while incorporating interactions with the surrounding ground. It ensures computational efficiency and differentiability, making it suitable for long-horizon optimization problems.

2. Formulation of a periodic OCP: This work introduces a novel OCP framework that jointly optimizes design parameters (e.g., installed capacities of renewable energy sources, thermal energy storage, and batteries) and optimal operational strategies to minimize total cost.

3. Reduction of computational complexity: To address the computational challenges associated with long-horizon optimization problems, an averaging method is applied to reduce computational complexity while preserving the accuracy needed for reliable solutions.

4. Economically feasible carbon-neutral energy systems: Through a case study of the Dietenbach district in Freiburg, Germany, the practical and economic feasibility of achieving carbon neutrality is demonstrated. The proposed system integrates renewable energy sources, and energy storage, providing a replicable model for sustainable energy systems at the district level.

## 2. System model

This chapter introduces a mathematical description of the system used in the optimal control problem. The proposed system is illustrated in Fig. 1, and the following subsections provide the details of each component.

### 2.1. Electric power balance

The electric power balance ensures that the total electricity input equals the total output at each time step. Mathematically, this balance can be expressed as:

$$P_{re}(t) - P_{load}(t) - P_{hp}(t) - P_b(t) + P_{grid}(t) = 0 \tag{2}$$

where

$$P_{re}(t) = s_{pv} P_{pv,0}(t) + s_{wind} P_{wind,0}(t) \geq 0 \tag{3}$$

represents the renewable power generation. Here, $P_{pv,0}(t)$ and $P_{wind,0}(t)$ are the known power profiles of photovoltaic (PV) and wind generation, considering the default installed capacities $C_{pv,0}$ and $C_{wind,0}$, respectively. These default installed capacities, scaled by respective design parameters $s_{pv}$ and $s_{wind}$, result in $C_{pv} = s_{pv}C_{pv,0}$, $C_{wind} = s_{wind}C_{wind,0}$. The battery power, $P_b(t)$, consists of charging $P_b^+(t)$ and discharging $P_b^-(t)$ variables, such that $P_b(t) = P_b^+(t) - P_b^-(t)$. Similarly, grid power exchange, $P_{grid}(t)$, is split into import $P_{grid}^+(t)$ and export $P_{grid}^-(t)$, expressed as $P_{grid}(t) = P_{grid}^+(t) - P_{grid}^-(t)$. The electricity demand profile $P_{load}(t)$ is treated as a known profile. All power variables are constrained to be non-negative:

$$0 \leq (P_b^+(t), P_b^-(t), P_{grid}^+(t), P_{grid}^-(t), P_{hp}(t)). \tag{4}$$

### 2.2. Battery model

The state of charge (SOC) of the battery, denoted $x^b(t) \in [0,1]$, evolves according to:

$$\dot{x}^b(t) = f^b\left(P_b^+(t), P_b^-(t); s_b\right) = \frac{1}{C_b}\left(P_b^+(t) \cdot \eta_{ch} - \frac{P_b^-(t)}{\eta_{dis}}\right). \tag{5}$$

Here, $C_b$ is the battery capacity in Wh, is scaled by $s_b$ based on a default capacity $C_{b,0}$. $\eta_{ch}$, $\eta_{dis}$ are the charging and discharging efficiencies,

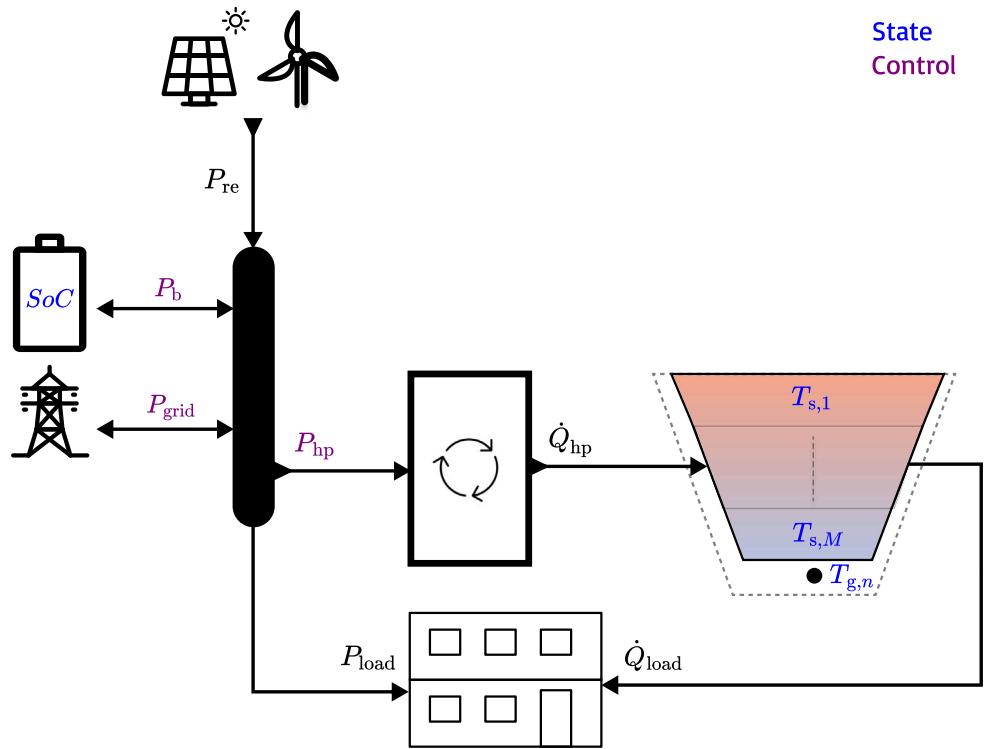

**Fig. 1.** System overview. State variables: storage temperature $T_{s,m}$, ground temperature $T_{g,n}$, and battery's state of charge; control inputs: power from the grid $P_{grid}$, battery $P_b$, and heat pump $P_{hp}$.

respectively, both assumed at 0.95. The SOC is constrained between 0 and 1:

$$0 \leq x^b(t) \leq 1. \qquad (6)$$

The charging and discharging power is limited to one-fourth of the total capacity per time step, corresponding to a 0.25C rate:

$$0 \leq P_b^+(t) \leq C_b/4\,\text{h}. \qquad (7a)$$

$$0 \leq P_b^-(t) \leq C_b/4\,\text{h}. \qquad (7b)$$

This model excludes detailed dynamics such as efficiency variations, self-discharge rates, and state-of-health degradation, as the primary focus of this paper is on seasonal thermal energy storage. More detailed battery models can be found in [25].

### 2.3. Heat pump model

The heat pump utilizes surplus renewable electricity to raise the temperature of thermal energy storage, rather than direct heating, for households. The coefficient of performance (COP) of the heat pump is:

$$COP = \eta_{Lorenz} \cdot \frac{T_{sink}}{T_{sink} - T_{source}} \qquad (8)$$

where $\eta_{Lorenz}$ is the Lorenz efficiency factor, set to 0.5 in this study to reflect industrial-scale heat pump performance [26], and $T_{source}$ corresponds to the ambient temperature $T_{amb}(t)$. The sink temperature is fixed at 86 °C, 1 K above the 85 °C storage limit, which keeps COP loss small and avoids numerical singularities in Eq. (12), where the mass flow term contains $(T_{sink} - T_s)$ in the denominator. This choice is further justified by the fact that industrial heat pumps with sink temperatures of 85 °C to 90 °C are common [27]. The heat transfer rate $\dot{Q}_{hp}$, is calculated as:

$$\dot{Q}_{hp}(t) = COP(t) \cdot P_{hp}(t). \qquad (9)$$

Heat pump's maximum capacity, defined as $C_{hp} = s_{hp} \cdot C_{hp,0}$, where $s_{hp}$ is a scaling parameter relative to the default capacity $C_{hp,0}$. The heat transfer rate is constrained by $C_{hp}$:

$$0 \leq \dot{Q}_{hp}(t) \leq C_{hp}. \qquad (10)$$

This study focuses on air-source heat pumps due to their increased efficiency in warmer weather, which aligns well with summer PV generation for charging seasonal storage. It should be noted that the corresponding model is simplified and omits operational details, such as startup dynamics, ramp-up times, or additional electrical power requirements. Including these effects would convert the formulated nonlinear programming problem below into a mixed-integer nonlinear programming problem, lengthening solution times by orders of magnitude [28]. Dynamic simulations show cycling losses of only 2% to 5% for cycle periods beyond 15 min [29]; with the 1h control step, the impact is negligible. These short-term transients are therefore neglected for the long-horizon optimization problem.

### 2.4. Thermal energy storage model and heating system

The model for the Pit Thermal Energy Storage (PTES) neglects the effects of groundwater flow (convective heat transfer in the ground) and focuses solely on heat conduction. The transient conduction process is simplified using a multi-node lumped-parameter approach, also known as the computational capacity resistance model (CaRM), which draws an analogy to electric circuits [30]. The concepts behind this model are illustrated in Fig. 2, which shows the heat transfer between the storage layers and the ground. Fig. 3 depicts the storage and ground behavior along with the thermal exchange between them.

#### 2.4.1. Geometry and basic assumptions

The geometry of the storage is a truncated pyramid, with square top and bottom faces of differing areas. Based on the district's heating demand, a default storage volume of 200 000 $m^3$ is adopted. This volume is identical to that of the world's largest operational PTES in Vojens,

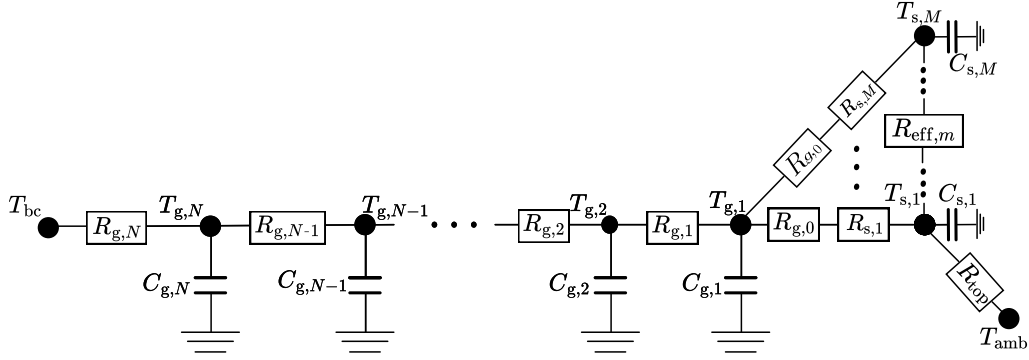

**Fig. 2.** Electric circuit analogy of heat transfer between heat storage and ground. Thermal capacitance $C$ represents individual storage or ground layers, temperature $T$ nodes are layer-average values, and thermal resistance $R$ connects adjacent nodes.

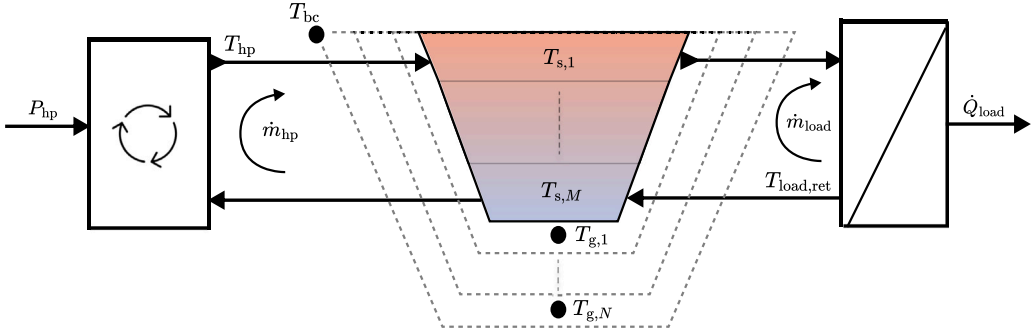

**Fig. 3.** Schematic of the multi-node pit storage model. The model consists of stratified storage layers and surrounding ground layers (dotted lines), with temperatures calculated at the center of each. During operation, heat is injected into or extracted from the top layer, while the return flow re-enters at the bottom to maintain stratification.

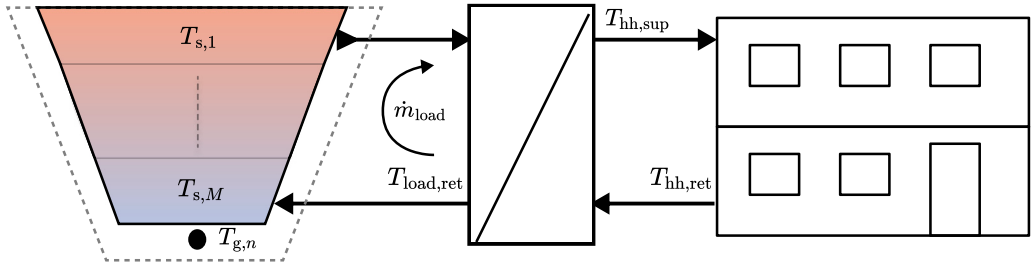

**Fig. 4.** Schematic representation of the heating supply and return flows in the system. The supply and return temperatures for the household heat demand, denoted as $T_{hh,sup}$ and $T_{hh,ret}$, are set to 40 °C and 20 °C, respectively.

Denmark, which serves a population and heat demand comparable to the case study (see Section 5.3) [31]. In the optimization, this reference volume is scaled by the factor $s_v$, allowing the solver to increase or decrease the size if that improves the objective. Initially, the side lengths of the top and bottom faces are 153.3 m and 73.2 m, respectively. The height is fixed at 15 m, a common choice for pit storage that reflects excavation cost constraints.

The underground geometry is approximated with planar walls, introducing inaccuracies at the edges; however, these are neglected in the model. For boundary conditions, a uniform ground temperature of 13.5 °C is assumed, based on the relatively stable temperatures observed at depths of 5 m to 10 m in the Freiburg region. More details of the geometry of the storage can be found in [14].

Fig. 4 illustrates the flow of heating supply to households and the corresponding return flow. The temperature of these flows fall within the typical range for low-temperature, fourth-generation district heating networks [32]. For simplicity, a constant temperature difference of 20 K between the supply and return flows is assumed. Thus, the return temperature of the heating load is consistently calculated as ($T_{load,ret} = T_{s,1} - 20$ K), assuming no heat exchange losses and a uniform heat transfer medium. This assumption simplifies the mass flow rate calculations, making $\dot{m}_{load}(t) = \frac{\dot{Q}_{load}(t)}{c_p \cdot 20\,\text{K}}$, which depends solely on $\dot{Q}_{load}(t)$.

### 2.4.2. Heat-transfer model and stratification

The storage temperature dynamics are governed by:

$$\dot{T}_{s,1} = \frac{1}{C_{s,1}} \left( c_p \cdot \dot{m}_{hp}(t)(T_{hp} - T_{s,1}) + c_p \cdot \dot{m}_{load}(t)(T_{s,2} - T_{s,1}) \right.$$
$$\left. + \frac{T_{s,2} - T_{s,1}}{R_{eff,1}} + \frac{T_{amb}(t) - T_{s,1}}{R_{top}} + \frac{T_{g,1} - T_{s,1}}{R_{g,0} + R_{s,1}} \right), \qquad (11a)$$

$$\dot{T}_{s,m} = \frac{1}{C_{s,m}} \left( c_p \cdot \dot{m}_{hp}(t)(T_{s,m-1} - T_{s,m}) + c_p \cdot \dot{m}_{load}(t)(T_{s,m+1} - T_{s,m}) \right.$$
$$\left. + \frac{T_{s,m-1} - 2T_{s,m} + T_{s,m+1}}{R_{eff,m}} + \frac{T_{g,1} - T_{s,m}}{R_{g,0} + R_{s,m}} \right), \quad \forall m = 2, \ldots, M-1, \qquad (11b)$$

**Table 1**
Parameters for the storage and underground model.

| Parameter | Value | Unit | Ref. |
|---|---|---|---|
| $\rho$ | 1000 | kg/m$^3$ | [12] |
| $c_p$ | 4200 | J/(kg K) | [12] |
| $\rho_g$ | 2000 | kg/m$^3$ | [12] |
| $c_{p,g}$ | 700 | J/(kg K) | [12] |
| $\lambda_g$ | 0.47 | W/(m K) | [12] |
| $\lambda_{eff}$ | 0.644 | W/(m K) | [33] |
| $U_{top}$ | 0.186 | W/(m$^2$ K) | [12] |
| $U$ | 90 | W/(m$^2$ K) | [13] |

$$\dot{T}_{s,M} = \frac{1}{C_{s,M}} \left( c_p \cdot \dot{m}_{hp}(t)(T_{s,M-1} - T_{s,M}) + c_p \cdot \dot{m}_{load}(t)(T_{load,ret} - T_{s,M}) \right.$$
$$\left. + \frac{T_{s,M-1} - T_{s,M}}{R_{eff,M}} + \frac{T_{g,1} - T_{s,M}}{R_{g,0} + R_{s,M}} \right), \tag{11c}$$

where

$$C_{s,m} = \rho c_p V_{s,m}, \quad \dot{m}_{hp}(t) = \frac{\dot{Q}_{hp}(t)}{c_p(T_{hp} - T_{s,M})}, \quad R_{s,m} = \frac{1}{U \cdot A_{surf,m}},$$

$$R_{eff,m} = \frac{z_m}{\lambda_{eff} \cdot A_{q,m}}, \quad R_{g,0} = M \cdot R_{g,0,original}, \quad R_{top} = \frac{1}{U_{top} \cdot A_{top}}. \tag{12}$$

Here, $C_{s,m}$ represents the thermal capacitance (in J/K) of each discretized section, calculated for a volume $V_{s,m}$, where all sections are divided into equal volumes. The stratified model accounts for forced convection due to $\dot{m}_{hp}(t)$ and $\dot{m}_{load}(t)$, and heat transfer between storage layers. The effective thermal resistance $R_{eff,m}$ is determined by the thermal conductivity $\lambda_{eff}$, the height $z_m$ of the $m$th layer, and its cross-sectional area $A_{q,m}$, assuming equal layer volumes. The thermal resistance $R_{s,m}$ accounts for the effective storage area of each layer, excluding the top layer, as it is insulated to minimize heat loss, while underground layers are uninsulated to utilize the ground's heat. Heat transfer between storage and ground layers is represented by $R_{g,0}$, adjusted for the number of layers $M$ due to their parallel connection. More details on the stratified model are provided in [33].

### 2.4.3. Ground-coupled heat transfer

The ground dynamics are modeled as:

$$\dot{T}_{g,1} = \frac{1}{C_{g,1}} \left( \frac{T_{s,1} - T_{g,1}}{R_{g,0} + R_{s,1}} + \cdots + \frac{T_{s,M} - T_{g,1}}{R_{g,0} + R_{s,M}} + \frac{T_{g,2} - T_{g,1}}{R_{g,1}} \right), \tag{13a}$$

$$\dot{T}_{g,n} = \frac{1}{C_{g,n}} \left( \frac{T_{g,n-1} - T_{g,n}}{R_{g,n-1}} + \frac{T_{g,n+1} - T_{g,n}}{R_{g,n}} \right), \forall n = 2, \ldots, N-1, \tag{13b}$$

$$\dot{T}_{g,N} = \frac{1}{C_{g,N}} \left( \frac{T_{g,N-1} - T_{g,N}}{R_{g,N-1}} + \frac{T_{bc} - T_{g,N}}{R_{g,N}} \right), \tag{13c}$$

where

$$C_{g,n} = \rho_g \cdot c_g \cdot V_n, \quad R_{g,n} = \frac{d_n}{\lambda_g \cdot A_{surf,n}}, \tag{14}$$

represent the thermal capacitance and resistance of the ground layers. Parameters for these dynamics are based on the PTES system in Dronninglund, Denmark [12], and are listed in Table 1. While temperature dynamics ($T$) and parameters such as thermal resistance ($R$) and capacitance ($C$) depend on the scaling factor ($s_x$), these dependencies are omitted in the equations for simplicity.

### 2.4.4. Discretization study

The numerical error of the lumped PTES model is governed by two geometric parameters: the number of radial ground layers $n$ and the outer-boundary distance $d$. To study these parameters with computational efficiency, the experiments employ a fully mixed storage model, which assumes a single temperature node for the entire water volume. One-year simulations with varying ($n, d$) were compared against a high-resolution benchmark ($n = 500$, $d = 100$ m) using the root mean squared

logarithmic error (RMSLE):

$$\text{RMSLE} = \sqrt{\frac{1}{N_f} \sum_{k=1}^{N_f} \left( \log_{10} T_s(n, d) - \log_{10} \bar{T}_s \right)^2} \tag{15}$$

which weights relative, not absolute, deviations and therefore prevents high-temperature periods from dominating the metric.

Fig. 5 shows that increasing the number of layers and distance lowers the RMSLE. With $n = 2$ and $d = 4$ m the error is already reduced to $0.5 \times 10^{-3}$, so finer ground meshes add only marginal benefit on an annual scale. Inside the storage, the vertical resolution is limited to four stratification layers ($m = 4$); Dahash et al. [12] demonstrated that further refinement mainly improves short-term thermocline detail while leaving yearly energy balances virtually unchanged. Accordingly, the grid configuration ($m, n, d$) = (4, 2, 4 m) is adopted for all optimization runs as a robust compromise between model fidelity and computational cost.

### 2.4.5. External validation

The first operating year of the Dronninglund PTES [31] is the only publicly available dataset for external validation. Its quality is constrained by (i) lack of bidirectional flow meters on the top and bottom diffusers—their SCADA signals were reconstructed from a mass- and energy-balance filter and still show a small imbalance, (ii) measurements covering only June–December 2014, and (iii) absence of near-wall soil temperatures, which prevents direct calibration of the ground model. The storage also underwent auxiliary flushing in the autumn, a process the model excludes.

To test the model's physical accuracy against this data, a detailed ten-layer configuration ($m = 10$, $n = 2$, $d = 4$ m) was first used to replicate the measured stratification. As shown in Fig. 6, this high-fidelity model reproduces the temperature dynamics with root mean squared errors of 2–3 °C for most layers.

For comparison, Fig. 7 shows the same validation scenario using the coarser four-layer model ($m = 4$) used in the main optimization study. While this simpler model still captures the main seasonal trends, it represents the thermal stratification with less detail than the ten-layer simulation.

The residual bias observed in both simulations is caused by uncertainties in the reconstructed flow data and unmodeled physical effects, such as corner diffusion. Although adding an interface resistance term could offer a marginal improvement, its influence is negligible for the four-layer model and is therefore omitted from the optimization study.

### 2.4.6. Summarized heat dynamics

The storage (including ground) temperatures form the state vector:

$$x^s = (T_{s,1}, \ldots, T_{s,M}, T_{g,1}, \ldots, T_{g,N})^\top \in \mathbb{R}^{N+M} \tag{16}$$

and the dynamics from (11) and (13) are given by:

$$\dot{x}^s = f^s(x^s, t; s_s). \tag{17}$$

The storage dynamics $f^s$ is time dependent and is influenced by $\dot{Q}_{hp}(t)$, $\dot{Q}_{load}(t)$, and $T_{amb}(t)$ with an affine dependence. Thus, the dynamics can be expressed in the form:

$$f^s(x^s, t; s_s) = f^s \left( x^s, \dot{Q}_{hp}(t), \dot{Q}_{load}(t), T_{amb}(t); s_s \right), \tag{18}$$

as described in Eqs. (11) to (13).

Storage temperatures are constrained between 10 °C and 85 °C, with the top storage layer maintaining a minimum temperature of 40 °C to meet the district heating network requirements. This can be expressed as:

$$40\,^\circ\text{C} \leq T_{s,1} \leq 85\,^\circ\text{C}, \tag{19a}$$

$$10\,^\circ\text{C} \leq T_{s,m} \leq 85\,^\circ\text{C}, \quad m = 2, \ldots, M. \tag{19b}$$

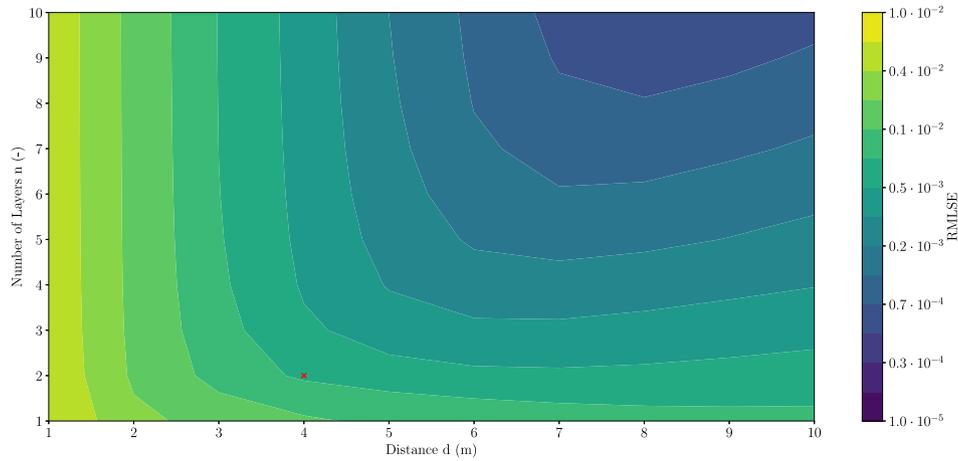

**Fig. 5.** RMSLE of the lumped PTES model as a function of ground-layer count $n$ and boundary distance $d$. The red cross marks the selected discretization $n = 2$, $d = 4$ m.

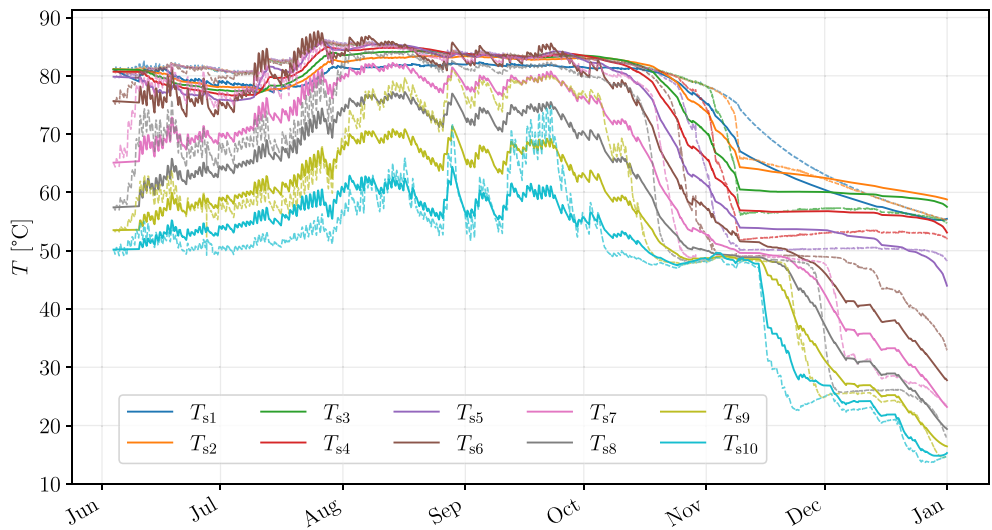

**Fig. 6.** Simulated (solid) and measured (dashed) layer temperatures for the Dronninglund PTES validation scenario (June–December 2014). The high-fidelity model shown here uses ten water layers and two ground layers.

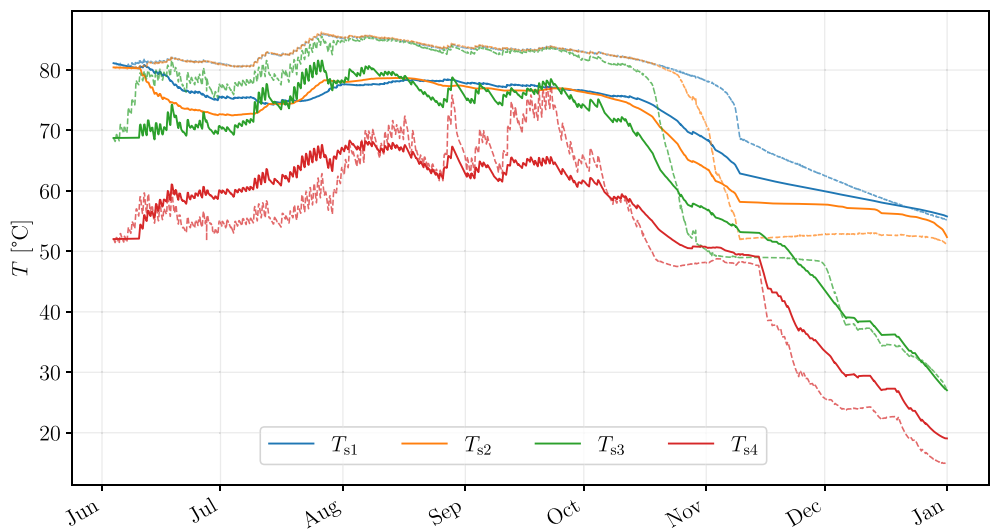

**Fig. 7.** Simulated (solid) and measured (dashed) layer temperatures for the four-layer storage model used in the main optimization study. This model is validated against the same Dronninglund PTES scenario shown in the previous figure.

**Table 2**

Investment parameters, and CAPEX/OPEX of each component[a].

|  | Inv. cost | Unit | Ref. | CAPEX | OPEX |
|---|---|---|---|---|---|
| PV | 1491 | $/kW$_\mathrm{p}$ | [34] | $I_\mathrm{pv} C_\mathrm{pv}$ | 0.01 CAPEX |
| Wind | 1569 | $/kW$_\mathrm{p}$ | [34] | $I_\mathrm{wind} C_\mathrm{wind}$ | 0.02 CAPEX |
| Battery | 476 | $/kWh | [34] | $I_\mathrm{b} C_\mathrm{b}$ | 0.02 CAPEX |
| Storage | 30 | €/m³ | [35] | $I_\mathrm{s} V_\mathrm{s}$ | 0.01 CAPEX |
| Heat pump | 651 | €/kW | [36] | $I_\mathrm{hp} C_\mathrm{hp}$ | 0.025 CAPEX |
| Total |  |  |  | CAPEX$_\mathrm{sum}$ | OPEX$_\mathrm{sum}$ |

[a] The exchange rate used is 1 $ corresponding to 0.92 € as the annual average for 2024.

## 3. Optimal control problem in continuous-time

This section presents a periodic optimization problem over a 30-year horizon aimed at minimizing the total cost, which includes both investment and operational costs of the proposed system. The objective is to optimize the design parameters and control strategy, assuming the system operates in a year-periodic steady state. For simplicity, the optimal control problem (OCP) is considered to have an infinite repeating cycle over time.

The OCP seeks to determine the fixed design parameters, which include the scaling factors of installed capacities of renewable energy sources, the heat pump, and the length of the storage area:

$$\theta = (s_\mathrm{pv}, s_\mathrm{wind}, s_\mathrm{b}, s_\mathrm{s}, s_\mathrm{hp}) \in \mathbb{R}^5. \tag{20}$$

Additionally, the OCP aims to optimize the trajectory of the control variables, which consist of the time-dependent power flows:

$$u(t) = \left( P_\mathrm{hp}(t), P_\mathrm{b}^+(t), P_\mathrm{b}^-(t), P_\mathrm{grid}^+(t), P_\mathrm{grid}^-(t) \right) \in \mathbb{R}^5, \quad \forall t \in [0, 365\mathrm{d}]. \tag{21}$$

The state vector $x(t)$ combines the temperatures and the battery's state of charge:

$$x(t) = \begin{bmatrix} x^\mathrm{s}(t) \\ x^\mathrm{b}(t) \end{bmatrix} \in \mathbb{R}^{M+N+1}, \quad \forall t \in [0, 365\mathrm{d}]. \tag{22}$$

The combined dynamics of the system are given by:

$$f(x(t), u(t), t; \theta) = \begin{bmatrix} f^\mathrm{s} \left( x^\mathrm{s}(t), \dot{Q}_\mathrm{hp}(t), \dot{Q}_\mathrm{load}(t), T_\mathrm{amb}(t); s_\mathrm{s} \right) \\ f^\mathrm{b}(P_\mathrm{b}^+(t), P_\mathrm{b}^-(t); s_\mathrm{b}) \end{bmatrix}. \tag{23}$$

Here, $\dot{Q}_\mathrm{hp}(t) = \dot{Q}_\mathrm{hp}(P_\mathrm{hp}(t), T_\mathrm{amb}(t), t)$ is a function of both the ambient temperature and the control $P_\mathrm{hp}(t)$, as defined in Eq. (9).

Fixed costs are divided into Capital Expenditure (CAPEX) and Operational Expenditure (OPEX), as summarized in Table 2. The annuity method is used to determine an annual payment for the CAPEX, ensuring that interest is covered each period and the total CAPEX is fully repaid by the end of the project duration ($n$ years), considering interest rate ($r$). The relationship between the CAPEX and the annual payment, referred to as the Annuity of the Investment (ANI), is expressed as:

$$\mathrm{CAPEX}_\mathrm{sum} = \mathrm{ANI} \cdot \sum_{t=1}^{n} \frac{1}{(1+r)^t}. \tag{24}$$

The geometric series can be expressed as:

$$\sum_{t=1}^{n} \frac{1}{(1+r)^t} = \frac{1 - (1+r)^{-n}}{r}. \tag{25}$$

Rearranging to solve for ANI one obtains:

$$\mathrm{ANI} = \mathrm{CAPEX}_\mathrm{sum} \cdot \frac{r(1+r)^n}{(1+r)^n - 1}. \tag{26}$$

The term $\frac{r(1+r)^n}{(1+r)^n - 1}$ is known as the Annuity Factor (ANF). For the cost calculation, $n = 30$ years and $r = 4\%$ are assumed. The fixed cost is then formulated as:

$$J_\mathrm{fix}(\theta) = \mathrm{ANI} + \mathrm{OPEX}_\mathrm{sum}. \tag{27}$$

The running cost is determined by electricity consumption from the grid and feed-in revenue:

$$J_\mathrm{run}(u(t)) = c_\mathrm{el,buy} P_\mathrm{grid}^+(t) - c_\mathrm{el,sell} P_\mathrm{grid}^-(t) \tag{28}$$

where $c_\mathrm{el,buy}$ and $c_\mathrm{el,sell}$ represent the electricity purchase price (30 ct/kWh) and feed-in tariff (1 ct/kWh), respectively. While feed-in tariffs are expected to decrease with the expansion of renewable energy deployment, the objective function is structured to promote grid independence by maintaining a high electricity purchase price. The continuous-time optimal control problem is summarized as:

$$\underset{x(\cdot), u(\cdot), \theta}{\mathrm{minimize}} \quad \int_0^{365\,\mathrm{d}} J_\mathrm{run}(u(t))\, dt \quad + \quad J_\mathrm{fix}(\theta) \tag{29a}$$

$$\text{subject to} \quad 0 = x(0) - x(365\,\mathrm{d}), \tag{29b}$$

$$\dot{x}(t) = f(x(t), u(t), t, \theta), \quad \forall t \in [0, 365\,\mathrm{d}], \tag{29c}$$

$$0 \geq h(u(t), t, \theta), \quad \forall t \in [0, 365\,\mathrm{d}], \tag{29d}$$

$$0 = g(u(t), t, \theta), \quad \forall t \in [0, 365\,\mathrm{d}], \tag{29e}$$

$$x_\mathrm{min} \leq x(t) \leq x_\mathrm{max}, \quad \forall t \in [0, 365\,\mathrm{d}], \tag{29f}$$

$$u_\mathrm{min} \leq u(t) \leq u_\mathrm{max}, \quad \forall t \in [0, 365\,\mathrm{d}], \tag{29g}$$

$$\theta_\mathrm{min} \leq \theta \leq \theta_\mathrm{max}. \tag{29h}$$

Here, (29b) is a periodicity constraint on the state trajectory. The inequality constraints (29d) consist of the parameter-dependent constraints (7) and (10), while the only equality constraint (29e) is given by the power balance (2). The state bounds (29f) are defined by Eqs. (19), and (6), while the control bounds (4). The parameter vectors are bounded in (29h), ranging from $\theta_\mathrm{min} = 0.1$ to $\theta_\mathrm{max} = 10$, to ensure better convergence and avoid numerical errors.

## 4. Discretization and nonlinear programming formulation

Solution strategies for optimal control problems (OCPs) are generally classified into indirect and direct methods. Indirect methods follow a "first optimize, then discretize" approach, require the derivation of first-order necessary conditions, and the solution of a boundary value problem [28]. In contrast, direct methods first discretize the state and control trajectories, transforming the problem into a finite-dimensional nonlinear programming (NLP) problem. This approach simplifies the treatment of state constraints and inequalities.

This work employs the direct multiple shooting method, where state trajectories are approximated by variables referred to as "shooting nodes" over a discretized time grid. Constraints derived from the system dynamics enforce consistency between these nodes. This approach is particularly well-suited for long-horizon optimization problems with complex dynamic constraints. Further details on transcribing OCPs using direct methods can be found in [28], while more information on solving the resulting NLPs is available in [37].

### 4.1. Discretization of optimal control problem

The numerical solution of the optimal control problem (29) begins by dividing the time horizon $t \in [0, 365\,\mathrm{d}]$ into $N_\mathrm{f}$ equidistant intervals, with a step size of $h_\mathrm{f} = 365\,\mathrm{d}/N_\mathrm{f}$, and grid node times $t_k = h_\mathrm{f} k, k = 0, \ldots, N_\mathrm{f} - 1$, such that:

$$0 = t_0 \leq t_1 \leq \cdots \leq t_k \leq \cdots \leq t_{N_\mathrm{f}-1} \leq t_{N_\mathrm{f}} = 365\,\mathrm{d}. \tag{30}$$

A typical selection is $N_f = 365 \cdot 24 = 8760$, corresponding to an hourly grid. This process establishes the "fine" grid; a corresponding "coarse" grid is introduced in Section 4.3.

The multiple shooting method discretizes the problem (29) by introducing shooting nodes $x_0, x_1, \ldots, x_{N_f} \in \mathbb{R}^{n_x}, x_k = (x_k^s, x_k^b)$ that approximate the state at the grid nodes as $x(t_k) \approx x_k$. The control variables $u_0, u_1, \ldots, u_{N_f-1}$ parameterize the control as piecewise constant: $u(t) = u_k \in \mathbb{R}^{n_u}, \forall t \in [t_k, t_{k+1})$. Additionally, the discretization treats the following quantities as constant over each grid interval:

$$
\begin{bmatrix} T_{amb}(t) \\ P_{pv}(t) \\ P_{wind}(t) \\ P_{load}(t) \\ \dot{Q}_{load}(t) \end{bmatrix} = \begin{bmatrix} T_{amb,k} \\ P_{pv,k} \\ P_{wind,k} \\ P_{load,k} \\ \dot{Q}_{load,k} \end{bmatrix}, \quad \forall t \in [t_k, t_{k+1}). \tag{31}
$$

The values of the parameters are provided from external data sources; further details can be found in Appendix A. The numerical implementation employs one of the simplest possible numerical integration methods for the multiple shooting interval: a single step of an implicit Euler integrator, which has proven to be the most efficient choice for this application. Note that more elaborate higher-order integrators exist and could prove to be more efficient in other applications. The system dynamics are then approximated as follows:

$$0 = x_{k+1}^s - x_k^s - h_f f^s(x_{k+1}^s, u_k, t_k, \theta), \qquad k = 0, 1, \ldots, N_f - 1 \tag{32}$$

$$0 = x_{k+1}^b - x_k^b - h_f f^b(u_k), \qquad k = 0, 1, \ldots, N_f - 1 \tag{33}$$

where Eq. (32) approximates the solution of the storage temperature dynamics (23) using an implicit Euler integration scheme of order 1. Eq. (33) for battery dynamics $f^b$ is exact under the applied piecewise-constant control parameterization. The integral of the cost function (28) for the piecewise constant control for the $k$th interval, is given by:

$$\int_{t_k}^{t_{k+1}} J_{run}(u(t)) = h_f J_{run}(u_k). \tag{34}$$

The transcribed continuous-time OCP (29) is then summarized as the following discretized NLP:

$$\underset{\substack{x_0, \ldots, x_{N_f}, \\ x_0^b, \ldots, x_{N_f}^b, \\ u_0, \ldots, u_{N_f-1}, \\ \theta}}{\text{minimize}} \sum_{k=0}^{N_f-1} h_f \cdot J_{run}(u_k) + J_{fix}(\theta) \tag{35a}$$

subject to
$$0 = x_0^s - x_{N_f}^s, \tag{35b}$$

$$x_{k+1}^s = x_k^s + h_f f^s(x_{k+1}^s, u_k, t_k, \theta), \quad k = 0, \ldots, N_f - 1, \tag{35c}$$

$$x_{min}^s \le x_k^s \le x_{max}^s, \qquad k = 0, \ldots, N_f, \tag{35d}$$

$$0 = x_0^b - x_{N_f}^b, \qquad k = 0, \ldots, N_f, \tag{35e}$$

$$x_{k+1}^b = x_k^b + h_f f^b(u_k), \qquad k = 0, \ldots, N_f - 1, \tag{35f}$$

$$x_{min}^b \le x_k^b \le x_{max}^b, \qquad k = 0, \ldots, N_f, \tag{35g}$$

$$0 \ge h(u_k, t_k, \theta), \qquad k = 0, \ldots, N_f - 1, \tag{35h}$$

$$0 = g(u_k, t_k, \theta), \qquad k = 0, \ldots, N_f - 1, \tag{35i}$$

$$u_{min} \le u_k \le u_{max}, \qquad k = 0, \ldots, N_f - 1, \tag{35j}$$

$$\theta_{min} \le \theta \le \theta_{max}. \tag{35k}$$

Here, the constraints are divided into groups: storage temperature states and dynamics (b–d), battery states and dynamics (e–g), control vector (h–j), and parameter constraints (k). The bounds on states and controls, as well as the (in)equalities $h$ and $g$, are defined as described in (29).

### 4.2. Reduced-size nonlinear programming formulation using the averaging method

The previously formulated NLP involves hourly time steps across an entire year, resulting in 8760 discretization points. Even for relatively simple models, this high dimensionality presents significant challenges for solvers. However, the temperature states $x(t)$ of seasonal thermal energy storage systems vary slowly, often over hours or days, even as control inputs fluctuate rapidly due to renewables or changes in heat demand. This property allows for a reduction in model complexity by averaging the storage temperature dynamics $f^s$ (17) over a daily time horizon $T$. The averaged dynamics for the storage state $x^s(t)$ over the next single day $t \in [t_k, t_k + T]$ are defined as:

$$\bar{f}_k^s(x^s(t)) = \frac{1}{T} \int_{t_k}^{t_k+T} f^s(x^s(t), \tau) \, d\tau \tag{36}$$

where, the averaging duration $T = 1\,\mathrm{d} = h_f K, K = 24$ is an integer multiple of the interval size $h_f$. This approach is justified using the Averaging Theorem, with further details provided in Appendix B.1.

The temperature dynamics $f^s$ depend on the time $t$ only via an affine dependence on the inputs $\dot{Q}_{hp}(t), \dot{Q}_{load}(t), T_{amb}(t)$, see Eq. (18). Exploiting the linearity of the averaging operator, the day-averaged dynamics (36) for $t \in [t_k, t_k + T]$ become:

$$\bar{f}_k^s(x^s(t)) = \frac{1}{T} \int_{t_k}^{t_k+T} f^s(x^s(t), \tau) \, d\tau \tag{37}$$

$$= \frac{1}{T} \int_{t_k}^{t_k+T} f^s\left(x^s(t), \dot{Q}_{hp}(\tau), \dot{Q}_{load}(\tau), T_{amb}(\tau)\right) \, d\tau \tag{38}$$

$$= f^s\left(x^s(t), \bar{\dot{Q}}_{hp,k}, \bar{\dot{Q}}_{load,k}, \bar{T}_{amb,k}\right) \tag{39}$$

where the values $\bar{\dot{Q}}_{hp,k}, \bar{\dot{Q}}_{load,k}, \bar{T}_{amb,k}$ represent the day-averaged heat transfer rates and ambient temperature, respectively. The detailed calculations for these values can be found in Appendix B.2.

### 4.3. Modifying the nonlinear programming problem with the averaged temperature dynamics

This methodology separates the simulation of the averaged temperature dynamics from the fast battery dynamics. To address the differing timescales between the slow averaged temperature dynamics and the fast battery dynamics, a second, coarse grid with $N_c \ll N_f$ intervals, indexed by $j = 0, \ldots, N_c$, is introduced. This coarse grid is related to the fine grid via $k = j \cdot K, K = 24$, and is used to integrate the slow, average temperature dynamics with large steps of size $h_c = K \cdot h_f = 24\,\mathrm{h}$. The averaging window size is chosen to match the step size of the slow integration scheme, which is set to one day ($T = h_c$). This process is illustrated in Fig. 8.

The constraints and integration equations for the averaged temperature states on the coarse grid are:

$$0 = x_0^s - x_{N_c}^s, \tag{40a}$$

$$x_{j+1}^s = x_j^s + h_c f^s\left(x_{j+1}^s, \bar{\dot{Q}}_{hp,k}, \bar{\dot{Q}}_{load,k}, \bar{T}_{amb,k}\right), \quad j = 0, \ldots, N_c - 1, \tag{40b}$$

$$x_{min}^s \le x_j^s \le x_{max}^s, \qquad j = 0, \ldots, N_c \tag{40c}$$

where $f^s$ is evaluated using an implicit Euler scheme with coarse step size $h_c$. The averaged inputs are calculated as in Eqs. (B.3), (B.4), and (B.7), with $k = jK$.

The modified NLP (35) is formulated as:

$$\underset{\substack{x_0, \ldots, x_{N_c}, \\ x_0^b, \ldots, x_{N_f}^b, \\ u_0, \ldots, u_{N_f-1}, \\ \theta}}{\text{minimize}} \quad (35a) \qquad \text{(Fine Grid Objective)} \tag{41a}$$

subject to
$$(40) \qquad \text{(Coarse Grid Constraints)}, \tag{41b}$$

$$(35e)–(35g) \qquad \text{(Fine Grid Constraints)}, \tag{41c}$$

$$(35h)–(35k) \qquad \text{(Control and Parameters Constraints)}. \tag{41d}$$

This formulation reduces the number of NLP variables from $N_f \cdot (M + N + 1 + n_u) + n_\theta$ to $N_c \cdot (M + N) + N_f \cdot (1 + n_u) + n_\theta$, achieving a reduction of approximately 50% in the example presented in Section 5.2. Here, $M$ and $N$ represent the number of storage and ground layers, respectively.

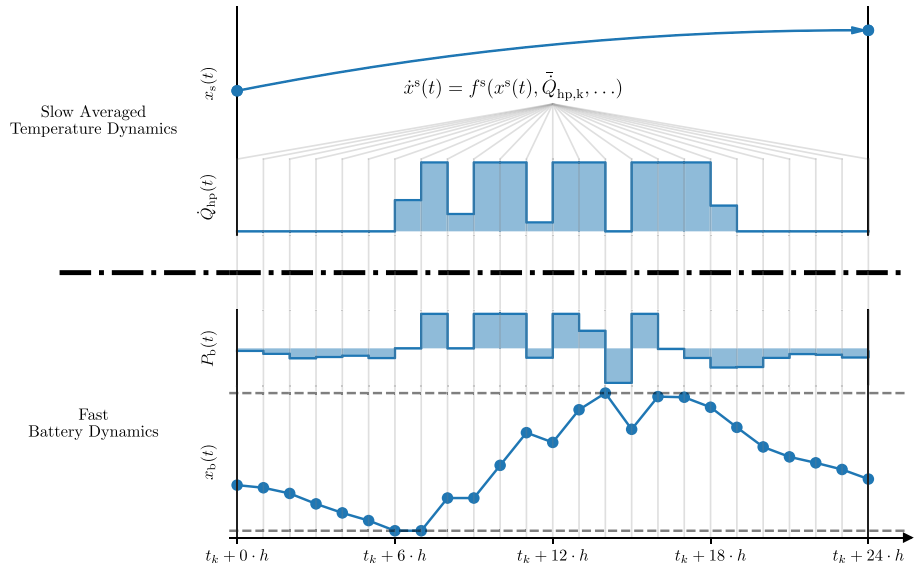

**Fig. 8.** The average dynamics of the states $x^s$ are integrated using coarse steps of size $h_c = 24$ h using the averaged inputs, while the fast battery dynamics are integrated over fine steps of $h_f = 1$ h to ensure constraint satisfaction.

This separation of numerical methods for different timescales is a well-established technique, applied in fields ranging from orbital-tether dynamics [38] and rotorcraft flight control [39] to flexible-spacecraft simulation [40].

*Numerical implementation and convergence checks.* The NLPs are formulated in CasADi [41] and solved using IPOPT [42] with the linear solver MA27 [43]. While the non-convex nature of these problems means a global optimum is not guaranteed, IPOPT converges reliably. The initial temperatures of the storage and ground are set to 30 °C and 13.5 °C, respectively. All controls are initialized to zero, and scaling parameters are set to one. The problem is scaled to ensure that all variables and their bounds are of order 1.

## 5. Results

In this section, the design and control optimization problem (29) is approximated by solving the nonlinear program (35). The optimization model is driven by hourly time series data for the Freiburg region for the year 2023. Air temperature data comes from the German Meteorological Service (DWD) [44]. Electricity demand is derived from the Energy Hub Design Optimization platform [4] and scaled to match the energy demand of the case study. Appendix A presents the detailed power profiles for PV, wind, and the heating load.

### 5.1. Optimal solution for a fixed electricity price

As a reference scenario, the optimization problem was solved for a constant electricity price of 0.3 €/kWh. Fig. 9 shows the optimal temperature profiles for the reference scenario. The results illustrate the seasonal storage performance and the thermal interaction between the storage and the surrounding ground.

The annual storage efficiency ($\eta_s$) is 94.3%, defined as the ratio of yearly useful heat delivered to the yearly heat charged:

$$\eta_s = \frac{\int_0^{365d} \dot{Q}_{out}(t)dt}{\int_0^{365d} \dot{Q}_{in}(t)dt}. \tag{42}$$

Here $\dot{Q}_{in}(t)$ is the heat supplied by the heat pump and $\dot{Q}_{out}$ is the heat delivered to the demand. The high efficiency is due to the large storage volume, which provides a low surface-to-volume ratio.

**Table 3**
Comparison of the numerical solution of full NLP (35) and averaged NLP (41) regarding complexity and the performance.

|                    | Original problem | Averaging method |
|--------------------|------------------|------------------|
| Total variables    | 105 130          | 54 762           |
| Time per iteration | 0.241 s          | 0.047 s          |
| Total solution time| 57.56 s          | 8.54 s           |

Additionally, a year-periodic steady state maintains higher surrounding ground temperatures, further boosting efficiency. For comparison, measured efficiencies of about 90% have been reported for the 60 000 m³ pit storage at Dronninglund [31], and numerical studies of large storage predict values above 90% [45].

### 5.2. Numerical study

By using Eq. (41), the optimal control problem can be solved with reduced time steps in the storage dynamics — 365 steps for one year — while maintaining the same objective function and number of control variables. Fig. 10 shows the optimal state trajectories of the averaged system compared to the original system. The differences are minimal, indicating that the averaged state dynamics closely match those of the original system.

This similarity is also evident in the optimal design parameters (maximum deviation: 1.15%), the optimal control trajectory (maximum deviation: 6.63% NRMSE[1] for $P_{hp}$) and the optimal state trajectory (maximum deviation: 3.08% NRMSE for $T_{g,1}$). Table 3 compares the computational complexity and performance between the solution for the full NLP (35) and the NLP (41) employing the averaging method.

### 5.3. Economic case study: Dietenbach

This case study assesses the economic feasibility of integrating PTES into a residential district currently under development in Freiburg, Germany.[2] The development will house approximately

---

[1] We use the normalized root mean squared error (NRMSE), defined as the RMSE normalized by the range of values: $\sqrt{1/N \sum (x_k - \hat{x}_k)^2}/(x_{max} - x_{min})$.

[2] Hydrogeological conditions in Freiburg, such as high groundwater flow, may affect PTES performance [46]. This limitation is discussed in Section 6.

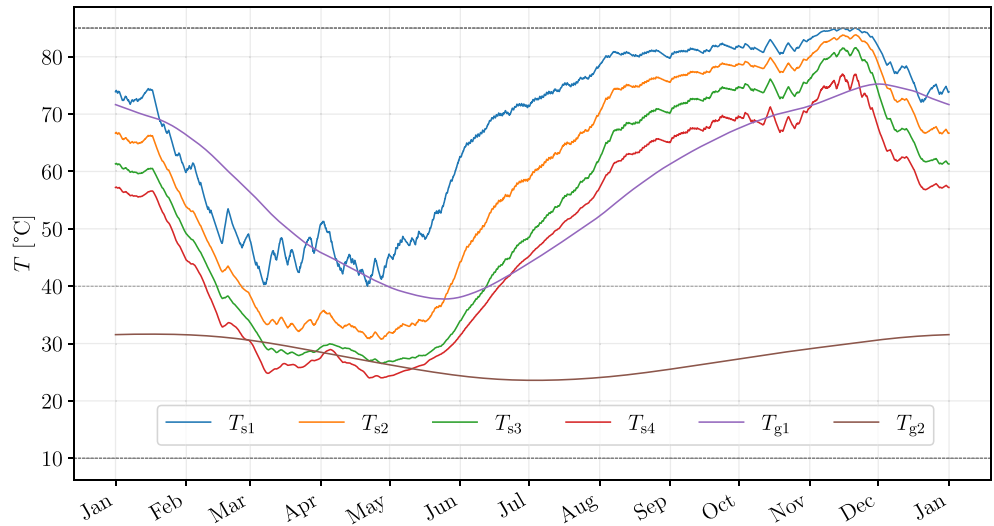

**Fig. 9.** Optimal temperature trajectories of the storage layers ($T_{s,1}$–$T_{s,4}$) and surrounding ground layers ($T_{g,1}$, $T_{g,2}$). The ground temperature oscillates in response to the storage, with the first layer ($T_{g,1}$) showing a stronger thermal connection than the second. The dashed line indicates the 40 °C cut-off temperature for heat discharge and the 85 °C maximum storage constraint.

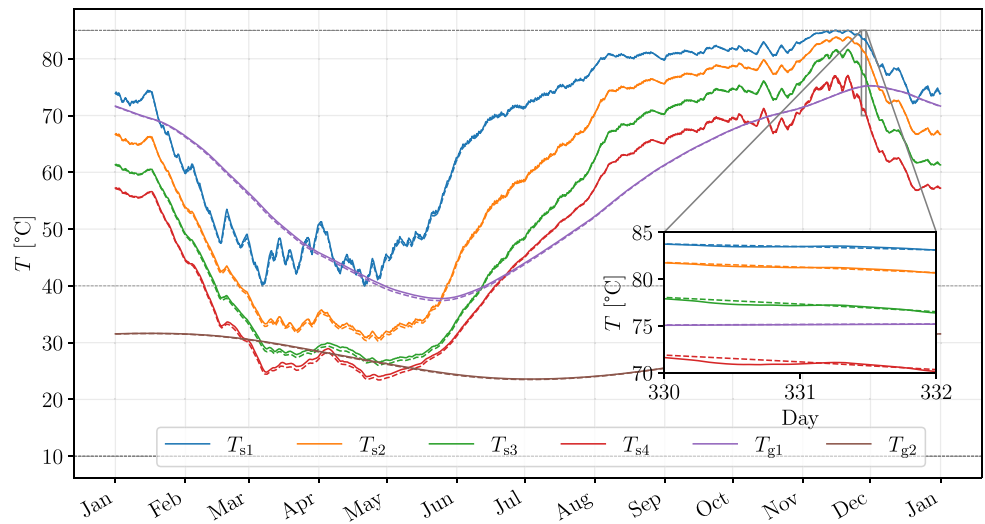

**Fig. 10.** Comparison of the optimal temperature trajectories of the full NLP (35) and the averaged NLP (41). The dotted lines represent the averaged state dynamics, while the solid lines represent the original dynamics.

16 000 residents across 6900 households, with an estimated annual heating demand of 51 GWh and electricity demand of 32 GWh [47]. Fig. 11 shows the weekly average demand profiles for heating and electricity.

### 5.3.1. Scenario comparison

Three system variants are now compared: (i) Full model — includes all components. (ii) No PTES — the full model without PTES. Heat is directly supplied by the heat pump at a supply temperature of 40 °C. (iii) No Wind — the full model without wind generation.

Table 4 compares the three variants. Normalized to Dietenbach's gross floor area of approximately 1.1 million m², the full model results in annual costs of 5.93 €/m². Although the district's buildings are well-insulated and thus require less heat, this cost is still considerably lower (approximately 73%) than the yearly German average of 22.1 €/m².[3]

A system without PTES, operating at the same 40 °C supply temperature, may benefit from a higher COP but fails to bridge the winter gap between renewable generation and heat demand effectively. As a result, excluding PTES leads to a 19.1% increase in total costs, clearly showing that integrating PTES is a cost-effective choice.

On the generation side, the complementary behavior of PV and wind, illustrated in Fig. 12, ensures a more consistent supply of renewable energy throughout the year: PV peaks in summer when heating demand is lowest, while wind output is strongest in winter, coinciding with peak heating needs. If wind generation is excluded, the system requires markedly higher capacities — PV (+101%), batteries (+23%), the heat pump (+102%), and PTES (+96%) — leading to a total cost increase of 45.9%, as illustrated in Table 4. This highlights the crucial role of wind generation in mitigating seasonal variability, even

---

[3] Household electricity use for appliances was 2893 kWh per dwelling in Germany [48]. Assuming an electricity price of 0.3 €/kWh and an average

dwelling size of 92.2 m² [49], this equals 9.4 €/m². Space heating demand was 131.4 kWh/m² [50]; for a heat pump with a seasonal COP of 3.1 [51], this corresponds to 12.7 €/m².

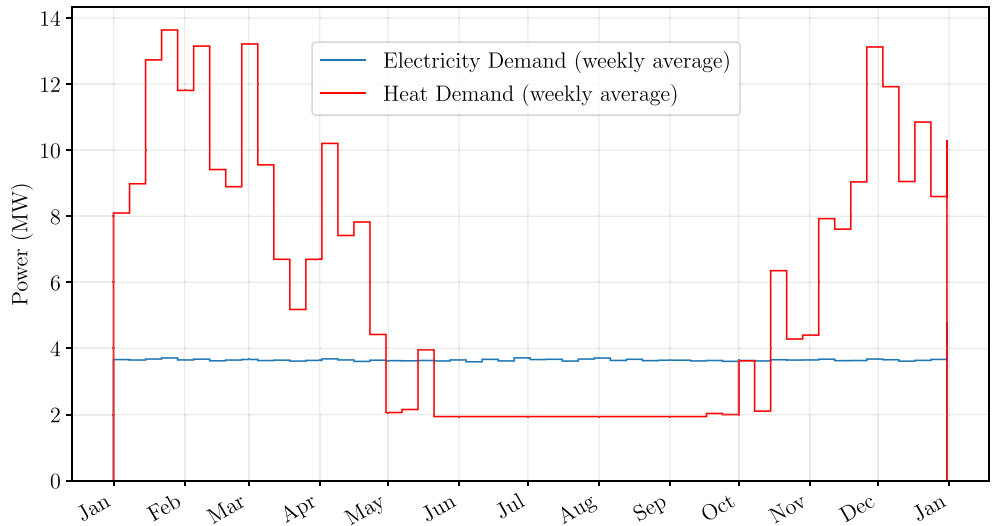

**Fig. 11.** Weekly averaged heating and electricity demand profiles for the Dietenbach district. The electricity demand shown here excludes the electricity consumed by heat pumps.

**Table 4**
Cost and capacity breakdown for three scenarios: (i) the full reference model, (ii) a system without PTES, and (iii) a system without wind power. All results are for a fixed electricity price of 0.3 €/kWh. All costs are presented as 30-year Net Present Values (NPVs) based on a 4% discount rate; operational costs (OPEX and Grid) are presented as their discounted present values.

| Category | Metric | (i) Full | (ii) No PTES | (iii) No wind | Unit |
|---|---|---|---|---|---|
| Grid | Import | 25.37 | 53.37 | 43.54 | M€ |
| | Export | −1.63 | −3.99 | −1.89 | M€ |
| Heatpump | Capacity | 13.86 | 20.03 | 27.96 | MW$_{th}$ |
| | CAPEX | 9.02 | 13.04 | 18.20 | M€ |
| | OPEX | 3.90 | 5.64 | 7.87 | M€ |
| PTES | Capacity | 261 851 | – | 512 564 | m³ |
| | CAPEX | 7.86 | – | 15.38 | M€ |
| | OPEX | 1.36 | – | 2.66 | M€ |
| PV | Capacity | 22.25 | 21.80 | 44.68 | MW$_p$ |
| | CAPEX | 30.51 | 29.89 | 61.25 | M€ |
| | OPEX | 5.28 | 5.17 | 10.59 | M€ |
| Wind | Capacity | 13.10 | 12.20 | – | MW$_p$ |
| | CAPEX | 18.91 | 17.61 | – | M€ |
| | OPEX | 6.54 | 6.09 | – | M€ |
| Battery | Capacity | 9.65 | 12.66 | 11.82 | MWh |
| | CAPEX | 4.23 | 5.54 | 5.18 | M€ |
| | OPEX | 1.46 | 1.92 | 1.79 | M€ |
| **Total** | CAPEX | 70.53 | 66.08 | 100.02 | M€ |
| | OPEX | 18.54 | 18.81 | 22.91 | M€ |
| | Total NPV | 112.79 | 134.27 | 164.57 | M€ |
| | Yearly cost | 5.93 | 7.06 | 8.65 | €/m² |

if its deployment encounters economic or political challenges in certain regions.

### 5.3.2. Impact of electricity price variations

The analysis now considers the impact of the electricity price variations on the system performance. Starting with an initial price of 0.1 €/kWh, $c_{el,buy}$ is increased by 0.1 €/kWh increments, up to 0.6 €/kWh. To compare scenarios, an "Autonomy" metric is defined as the ratio of renewable energy generation to total electricity demand:

$$\text{Autonomy} := \frac{E_{tot} - E_{grid,buy}}{E_{tot}} \tag{43}$$

where,

$$E_{tot} = \int_0^{365d} (P_{hp}(t) + P_{load}(t)) \, dt \tag{44}$$

$$E_{grid,buy} = \int_0^{365d} P_{grid}^+(t) \, dt \tag{45}$$

As shown in Fig. 13, increasing grid prices leads to higher total system costs and improved autonomy levels. However, full autonomy is not reached because renewable generation during winter months is insufficient, even with significant investments in battery and renewable energy capacity driven by high electricity prices. Consistent with this finding, scenarios without wind generation or PTES exhibit significantly reduced autonomy and increased costs, as these cases cannot adequately address seasonal variability in renewable energy generation.

## 6. Conclusions

This work developed a mathematical model of a renewable energy system, including seasonal thermal energy storage (STES), a heat pump, a battery, and the electricity grid. The underground storage system was represented using a one-dimensional, multi-node lumped-parameter scheme with uniform boundary conditions, capturing the mass flow and stratified temperature dynamics. An optimization study based on these models yielded the following key findings:

1. **Real-world Case Study**: For the Dietenbach district in Freiburg, Germany, the proposed system reached 91% self-sufficiency and cut household energy cost by approximately 73% relative to the German average.
2. **Efficient Optimal Control**: The averaging method reduced computation time by 80.5% with minimal deviation from the detailed model.

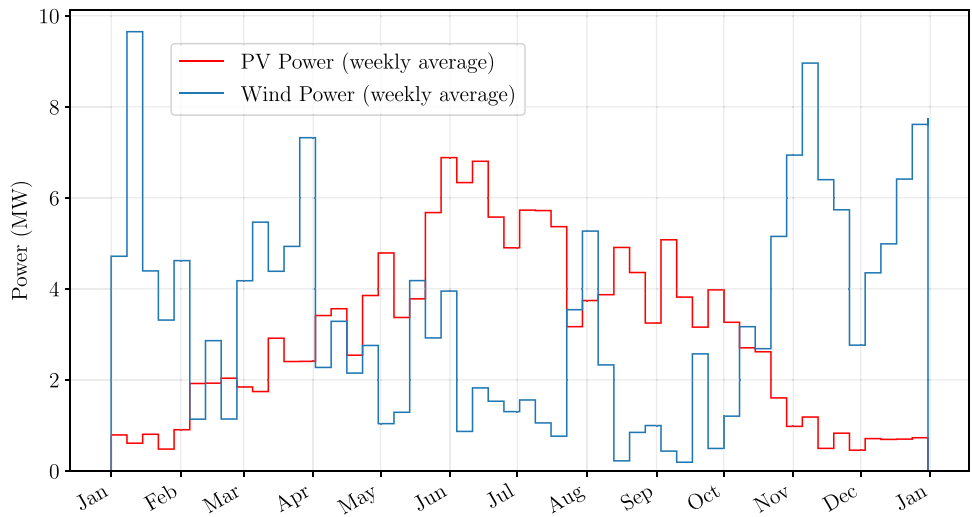

**Fig. 12.** Weekly averaged generation profiles for PV and wind, showing their opposing seasonal trends.

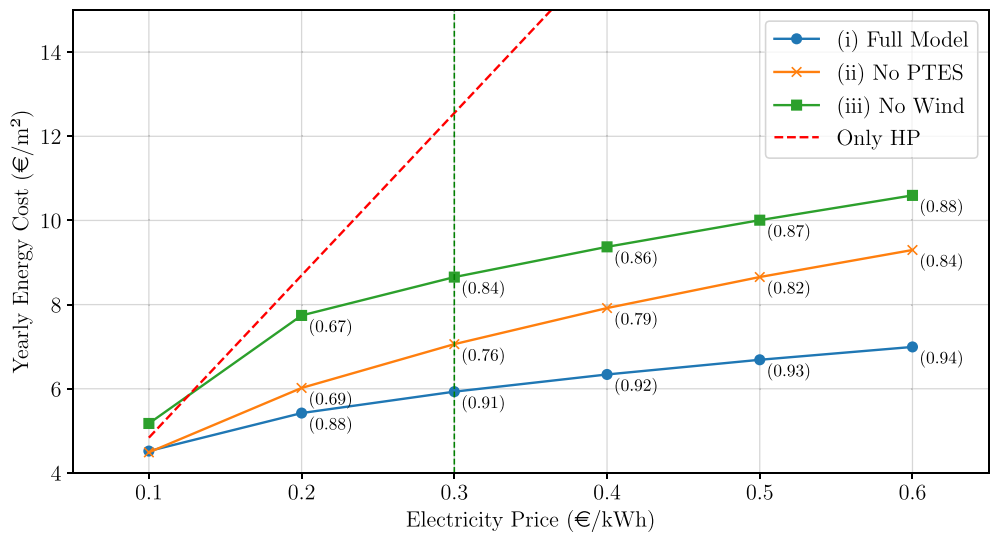

**Fig. 13.** Cost comparison based on electricity price variations, with corresponding autonomy levels shown in brackets. The only HP scenario consists of a heat pump with no local renewable generation, fully dependent on the grid to meet its energy demand.

The reported figures represent an optimistic upper bound, as several simplifying assumptions were introduced to ensure computational tractability. In the component models, heat-pump start-up transients, ramp-rate limits, and minimum on–off periods were excluded, resulting in a best-case dispatch scenario. Additionally, part-load efficiencies were neglected, and component capacities were treated as continuous rather than discrete variables. For the STES model, the discretization of both the storage and ground layers employed a limited node count to balance fidelity with computation time. Geological heterogeneity, the effects of groundwater flow, and variations in thermal diffusivity were also omitted. In particular, while conduction typically dominates

STES heat losses, recent work shows that groundwater can considerably increase them under certain hydrogeological conditions [52]. For the Freiburg case study, this simplification implies that the economic results should be regarded as an optimistic baseline. The dynamics of the district heating network itself were simplified, which is a limitation for detailed network planning but acceptable for assessing system-wide scalability. Finally, the analysis was performed under a deterministic framework with fixed meteorological conditions, electricity prices, and load profiles, meaning that input uncertainty was not captured.

Future research could enhance the model by incorporating uncertainty in weather and load data, as well as integrating more detailed dynamics for the heat pump and the district heating network. Detailed investigations of site-specific geological conditions are essential to accurately assess STES performance in practice. Furthermore, advanced control strategies could be explored, such as adopting a flexible use of the storage layers, using the heat-pump sink temperature as a control variable, and integrating domestic hot water boosts to improve efficiency and cost-effectiveness.

Overall, this study broadens the understanding of STES beyond traditional solar thermal applications, demonstrates its feasibility in real-world settings, and supports the growing imperative for electrification and decarbonization of the heating sector. It also advances research on optimal control for renewable energy systems and highlights the practical application of the averaging method for long-horizon optimization.

## CRediT authorship contribution statement

**Wonsun Song:** Writing – review & editing, Writing – original draft, Visualization, Validation, Software, Methodology, Formal analysis, Data curation, Conceptualization. **Jakob Harzer:** Writing – review & editing, Writing – original draft, Visualization, Validation, Software, Methodology, Conceptualization. **Christopher Jung:** Data curation. **Leon Sander:** Data curation. **Moritz Diehl:** Supervision, Conceptualization.

## Open source implementation

The mini-toolbox for the code is available at https://github.com/JakobHarz/NOSTES; it includes all input data and the scripts required to reproduce the optimization and validation results presented in this work.

## Declaration of competing interest

The authors declare that they have no known competing financial interests or personal relationships that could have appeared to influence the work reported in this paper.

## Acknowledgments

This research was supported by the German Research Foundation (DFG) via project 525018088, and by the Federal Ministry for Economic Affairs and Energy of Germany (BMWK) via 03EN3054B. The authors thank Manuel Kollmar, Armin Nurkanović, Dirk Schindler, and Arne Groß for their help, guidance, and fruitful discussions.

## Appendix A. Modeling details

### A.1. Photovoltaics model

The hourly power output of a solar PV module in Freiburg is determined using hourly-averaged global horizontal irradiance (GHI) and air temperature ($T_{amb}$) data from the German Meteorological Service station [44]. The power output ($P_{out,pv}$) is calculated based on [53]:

$$P_{out,pv} = GSI \cdot \eta_{ref} \left[1 - \beta_{ref}(T_{cell} - T_{ref})\right] \cdot PR \quad (A.1)$$

where the irradiance on the module's surface ($GSI$) is given in W/m². For an assumed 30° south-facing tilt, it is estimated by dividing the GHI data by cos(30°) [54]. The module's cell temperature ($T_{cell}$), in °C, is estimated empirically:

$$T_{cell} = c_1 + c_2 T_{amb} + c_3 GSI \quad (A.2)$$

where the coefficients are $c_1 = -3.75\,°C$, $c_2 = 1.14$, and $c_3 = 0.0175\,°C\,m^2/W$.

The specific parameters are for a modern Longi Hi-Mo X6 Scientist LR5-72 HTH 590-600M [55], which has an electrical efficiency ($\eta_{ref}$) of 23.2%, a temperature coefficient ($\beta_{ref}$) of −0.290%/°C, and a reference temperature ($T_{ref}$) of 25 °C under standard testing condition. A performance ratio ($PR$) of 0.9 accounts for system losses such as wiring and soiling, based on literature values [56].

### A.2. Wind model

The wind turbine power output is estimated using hourly wind speed data simulated at a height of 166 m from the Wind speed Complementary Model (WiCoMo) [57] and the power output of a 5.6 MW wind turbine. Although the model accounts for integer numbers of turbines, the power output is scaled by $s_{wind}$ to represent fractional installed capacities for the optimization study.

### A.3. Heating load model

The heating demand is divided into space heating (SH) and domestic hot water (DHW) heating. The space heating calculation uses simplified model that incorporates the effective temperature, $T_{eff}$, calculated as the 1-day moving average of $T_{amb}$:

$$T_{eff}(t) = \frac{1}{1d} \int_{t-1d}^{t} T_{amb}(\tau)\,d\tau \quad (A.3)$$

Heating is supplied when $T_{eff}(t)$ falls below the threshold temperature, $T_{border} = 12\,°C$. The heating demand depends on the temperature difference between the room temperature $T_{room} = 20\,°C$ and $T_{eff}(t)$. The space heating demand is scaled by the factor $s_{sh}$, resulting in a total space heating demand of 34 GWh. This relationship is expressed as:

$$\dot{Q}_{sh}(T_{eff}(t)) = \begin{cases} 0, & \text{if } T_{eff}(t) \geq T_{border}, \\ s_{sh} \cdot (T_{room} - T_{eff}(t)), & \text{if } T_{eff}(t) < T_{border}. \end{cases} \quad (A.4)$$

The domestic hot water demand, $\dot{Q}_{dhw}$, is modeled with daily peak usage in the morning and evening. These peaks are represented by Gaussian functions:

$$\dot{Q}_{dhw}(t) = s_{dhw} \cdot (\dot{Q}_{morning}(t) + \dot{Q}_{evening}(t)), \quad (A.5)$$

where $\dot{Q}_{morning}(t)$ and $\dot{Q}_{evening}(t)$ are defined as:

$$\dot{Q}(t) = \exp\left(-\frac{(t-\mu)^2}{2\sigma^2}\right), \quad (A.6)$$

with $\mu$ and $\sigma$ representing the mean and standard deviation of the peak demand, respectively. The demand is scaled by $s_{dhw}$, which is 17 GWh in our case study.

Finally, the total heating demand, $\dot{Q}_{load}(t)$, is the sum of the space heating and domestic hot water demands, 51 GWh:

$$\dot{Q}_{load}(t) = \dot{Q}_{sh}(t) + \dot{Q}_{dhw}(t). \quad (A.7)$$

## Appendix B. Averaging method

### B.1. Theory of the averaging method

The use of average dynamics is motivated by the following theorem, adapted from [58] (Lemma 2.1.8),

**Lemma 1** (*Averaging Theorem*). *Consider the solution $x(t)$ of the initial value problem*

$$\dot{x}(t) = \epsilon f(x(t), t), \qquad x(0) = a \quad (B.1)$$

*where the dynamics $f^1 : \mathbb{R}^{n_x} \times \mathbb{R}$ are Lipschitz continuous in $(x,t)$ and $T$-periodic in $t$, and $0 < \epsilon \ll 1$. Then, the solution $y(t)$ of the averaged dynamics*

$$\dot{y}(t) = \epsilon \bar{f}(y(t)) = \epsilon \frac{1}{T} \int_0^T f(y(t), \tau)\,d\tau \quad (B.2)$$

*with initial value $y(0) = a$ approximates the solution of the original dynamics as $y(t) = x(t) + \mathcal{O}(\epsilon)$ on the timescale $1/\epsilon$.*

Sanders et al. [58] provide a proof and further details. Some remarks are in order. First of all, the temperature dynamics $f^s$ (Eq. (17)) are of the particular form $\epsilon f$ where the parameter $\epsilon$ is equivalent to the inverse of the thermal capacity of the smallest discretization volume. Secondly, although the dynamics $f^s(x, t)$ are not 1-periodic, the time-dependence is treated as locally periodic because the averaging is only computed for a particular day: $f^s_k(x^s, t) = f^s(x^s, t_k + (t \mod T))$. Third, since the controls and time series for the ambient temperature and heat transfer rates are applied as piecewise constant over the intervals, the dynamics $f_s$ are not Lipschitz continuous in $t$. Nevertheless, the averaging theorem still holds; a proof of the averaging theorem without this restriction can be found in [59].

### B.2. Calculation details for the averaging method

The day-averaged values of $\dot{Q}_{load}(t)$ and $T_{amb}(t)$ are calculated as single sums, due to the piecewise constant parameterization and the averaging horizon $T$ being integer multiples of the step size $h_f$, given by:

$$\bar{Q}_{load,k} := \frac{1}{T} \int_{t_k}^{t_k+T} \dot{Q}_{load}(t)\, dt = \frac{1}{K} \sum_{i=k}^{i=k+K-1} \dot{Q}_{load,i} \tag{B.3}$$

$$\bar{T}_{amb,k} := \frac{1}{T} \int_{t_k}^{t_k+T} T_{amb}(t)\, dt = \frac{1}{K} \sum_{i=k}^{i=k+K-1} T_{amb,i}. \tag{B.4}$$

For $\dot{Q}_{hp}(t)$, averaging requires careful consideration of its non-linear dependence on $T_{amb}$. It is given by:

$$\bar{Q}_{hp,k} := \frac{1}{T} \int_0^T \dot{Q}_{hp}(t)\, dt \tag{B.5}$$

$$= \frac{1}{T} \int_0^T \eta_{Lorenz} \cdot \frac{T_{hp}}{T_{hp} - T_{amb}(t)} \cdot P_{hp}(t)\, dt \tag{B.6}$$

$$= \frac{1}{K} \sum_{i=k}^{i=k+K-1} \eta_{Lorenz} \cdot \frac{T_{hp}}{T_{hp} - T_{amb,i}} \cdot P_{hp,i} \tag{B.7}$$

Here, the term $\eta_{Lorenz} \cdot \frac{T_{hp}}{T_{hp}-T_{amb,i}}$ represents time-dependent but fixed problem data, while $P_{hp,i}$ is optimization variable.

### Data availability

The mini-toolbox for the code is available at https://github.com/JakobHarz/NOSTES; it includes all input data and the scripts required to reproduce the optimization and validation results.